\documentclass[12pt]{amsart}
\usepackage{amsmath}
\usepackage{amsfonts,amssymb,amsthm, txfonts, pxfonts}
\def\struckint{\mathop{%
\def\mathpalette##1##2{\mathchoice{##1\displaystyle##2}%
  {##1\textstyle##2}{##1\scriptstyle##2}{##1\scriptscriptstyle##2}}%
\mathpalette
{\vbox\bgroup\baselineskip0pt\lineskiplimit-1000pt\lineskip-1000pt
\halign\bgroup\hfill$}
{##$\hfill\cr{\intop}\cr\diagup\cr\egroup\egroup}%
}\limits}

\newcommand{\essn}{{\mathbb{S}^n}}

\newtheorem{theorem}{Theorem}
\newtheorem{lemma}[theorem]{Lemma}
\newtheorem{corollary}[theorem]{Corollary}

\newtheorem{definition}[theorem]{Definition}

\newtheorem{theorem-definition}[theorem]{Theorem-Definition}

\theoremstyle{remark}
\newtheorem{observation}[theorem]{Observation}
\newtheorem{remark}[theorem]{Remark}

\newtheorem{question}[theorem]{Question}

\newcommand{\cx}{\mathbb{C}}

\newcommand{\reals}{\mathbb{R}}

\DeclareMathOperator{\diag}{Diag}
\DeclareMathOperator{\vol}{Vol}
\DeclareMathOperator{\tr}{tr}
\DeclareMathOperator{\acosh}{arccosh}

\begin{document}


\title{On some mean matrix inequalites of dynamical interest}

\author{Igor Rivin}

\address{Department of Mathematics, Temple University, Philadelphia}

\email{rivin@math.temple.edu}

\thanks{The author is supported by the NSF DMS. The author would like to thank
  Amie Wilkinson for bringing the question to his attention, and Omar Hijab
  for suggesting the current simple proof of Theorem \ref{redineq}}

\date{\today}



\begin{abstract}
Let $A\in SL(n, \mathbb{R}).$ We show that 
for all $n > 2$ there exist dimensional strictly positive constants $C_n$
such that 
\[
\fint_{O_n} \log \rho(A X) d X \geq C_n \log \|A\|,
\]
where $\|A\|$ denotes the operator norm of $A$ (which equals the largest
singular value of $A$), $\rho$ denotes the spectral radius, and the integral
is with respect to the Haar measure on $O_n.$ The same result (with
essentially the same proof) holds for the unitary group $U_n$ in place of the
orthogonal group. The result does not hold in dimension $2.$
This answers questions asked in \cite{bpsw01,lssw03,ds01}. We also discuss
what happens when the integral above is taken with respect to measure other
than the Haar measure.

We also give a simple proof that 
\[
\int_{\mathbb{S}^n} \log \|A u\| d\sigma_n \geq 0,
\]
with equality if and only if $A^t A = I_n,$ ($d\sigma_n$ is the standard
measure on $\essn.$)

\end{abstract}

\maketitle

\section*{Introduction}

A.~Wilkinson brought the following questions to the author's attention:
\begin{question}
\label{wilk}
Let $A \in SL(n, \mathbb{R}).$ Consider the coset $\mathcal{A} = A O(n).$
Then is it true that 
\[
\int_{\mathcal A} d\mathfrak{o} \int_\essn \log \|B u\| d d\sigma_n\neq 0,
\] 
whenever $A$ is not orthogonal? (In the above
integral, $d\mathfrak{o}$ denotes the Haar measure on the orthogonal group,
normalized to be a probability measure; $d\sigma_n$ is the standard measure on
$\essn.$)
\end{question}
A version of this question is addressed in the paper \cite{lssw03}. In this
note we give a very simple proof that the answer is in the affirmative
(Theorem \ref{mainineq}).

\begin{question}
\label{wilk2}
Let $A \in SL(n, \reals),$ then
\[
\int_{O_n} \log \rho(A X) d X \geq C_n \int_{\mathcal A} d\mathfrak{o} \int_\essn
\log \|B u\| d d\sigma_n\neq 0,
\]
for some positive dimensional constant $C_n.$
\end{question}
The above is Question 6.6 in \cite{bpsw01}. The question is answered in
dimension $2$ in the paper \cite{ab01} -- in fact, the authors show that in
that case one can take $C_n = 1,$ and replace $\geq$ by $=.$

In the papers \cite{bpsw01,lssw03} the question is also posed as to what
happens when the measure on the orthogonal group is \emph{not} the Haar
measure (actually, the question is asked in the slightly different context.)
Obviously the result is unchanged if a measure $\mu$ is in the same measure
class as the Haar measure, but we note that in \emph{two} dimensions the sign of the
inequality in Question \ref{wilk2} can be reversed whenever $\mu$ is supported
on a proper subset of $O(2) = S^1$ -- see Section \ref{dim2} for details.

The same question with the orthogonal group replaced by the unitary group is
answered in the paper \cite{ds01} (in arbitrary dimension). The dimensional
constant is again equal to $1,$ though the inequality cannot be replaced by an
equality. In \cite{ds01}, Dedieu and Shub also conjecture a positive answer to
Question \ref{wilk2}. The same conjecture is made in \cite{lssw03}.

In this note we give a simple proof of 
\begin{theorem}
\label{mainineq}
Let $A\in SL(n, \mathbb{R}).$ Then
\[
\int_{\mathbb{S}^n} \log \|A u\| d\sigma n \geq 0,
\]
with equality if and only if $A^t A = I_n,$ (and $d \sigma_n$ denotes the
standard measure on the unit sphere).
\end{theorem}
and, as one of the consequences, we answer Question \ref{wilk} in the
affirmative.

We also show 
\begin{theorem}
\label{moremainineq}
For all $n > 2$ there exist strictly positive constants $C_n$
such that 
\[
\fint_{O_n} \log \rho(A X) d X \geq C_n \log \|A\|,
\]
where $\|A\|$ denotes the operator norm of $A$ (which equals the largest singular
value of $A$), $\rho$ denotes the spectral radius, and the integral is with
respect to the Haar measure on $O_n.$ The same result (with essentially the
same proof) holds for the unitary group $U_n$ in place of the orthogonal
group.
\end{theorem}
This result easily implies a positive answer to Question \ref{wilk2} in
dimensions bigger than $2.$ Theorem \ref{moremainineq} is \emph{false} in
dimension $2.$ Some discussion of why dimension $2$ is different is given in
Section \ref{dim2}.

It should also be noted that Theorem \ref{moremainineq} immediately implies
the following 
\begin{corollary}
\label{genmu}
Let $\mu$ be an orthogonally invariant measure on $SL(n, \reals),$ such that 
$\log \|x\| \in L^1(\mu),$ and $n \geq 3.$
Then 
\[
\int_{SL(n, \reals} \log \rho( X) d \mu \geq C_n \int_{SL(n, \reals)}\log
\|A\| d\mu. 
\]
The same statement holds with $\cx$ replacing $\reals$ and ``unitarily
invariant'' replacing ``orthogonally invariant''.
\end{corollary}
\begin{proof}
Integrate both sides of the inequality in Theorem \ref{moremainineq}.
\end{proof}


The plan of the paper is as follows: In Section \ref{generalities} we
summarize some necessary facts of linear algebra. In Section \ref{mainproof}
we give a proof of Theorem \ref{mainineq}. In Section \ref{lbound} we outline
our approach to the proof of Theorem \ref{moremainineq}. The actual proof
falls naturally into two parts -- one part covers the case where $A^*A$ is far
from $I_n$ and is the content of Section \ref{away} and the other part the
case where it is close to the identity, and is the content of Section
\ref{hiho}. This last section requires some simple results on perturbation
theory, and these are covered in Section \ref{perturb}.

It should be remarked that no effort is made to estimate the constant $C_n$
which appears in the statement in the Theorem \ref{moremainineq} -- in order
to have a sharp estimate, the methods of the current paper will need to be
combined with detailed analysis of the structure of the orthogonal groups, and
this will be the subject of another paper.

\section{Generalities on matrices}
\label{generalities}
First,
\begin{definition}
Let $A \in M_{m\times n}.$ Then the nonnegative square roots of the
eigenvalues of the (positive semi-definite) matrix $A^tA$ are called the
\emph{singular values} of $A.$  We will denote individual singular values by
$\sigma_1\geq \dots \geq  \sigma_m,$ and we shall denote the ordered $m$-tuple
$(\sigma_1, \dotsc, \sigma_m)$ by $\Sigma(A).$
\end{definition}

In this note we will only work with square matrices (although the results go
through without change for non-square matrices), and our use of singular
values will be localized to the following two lemmas:

\begin{lemma}
\label{invar}
Let $A\in M_{n\times n},$ and let  $P \in O(n).$ Then 
\[
\Sigma(A) = \Sigma(A P) = \Sigma(P A).
\]
\end{lemma}

\begin{proof} Immediate. \end{proof}

\begin{lemma}
\label{singinv}
Let $A \in M_{n\times n},$ and let $f \mathbb{R} \rightarrow \mathbb{R}$ be
such that 
\[
\mathfrak{I} = \int_\essn f\left(\langle A u, A u \rangle\right) d\sigma_n < \infty.
\]
Then 
\[
\mathfrak{I} = \int_\essn f\left(\sum_{i=1}^n \sigma_i^2(A) u_i^2\right)d \sigma_n.
\]
\end{lemma}
\begin{proof}
It is sufficient to note that 
\[
\langle A u, A u\rangle = u^t A^t A u = (Pu)^t \diag\left(\Sigma^2(A)\right)
Pu,
\]
where $P$ is an orthogonal matrix (hence an isometry of $\essn$).
\end{proof}
\begin{corollary}
\[
\int_\mathcal{A} d\mathfrak{o} \int_\essn \log \|B u\| d \sigma_n = 
\dfrac{1}{2} \int_\essn \log \left(\sum_{i=1}^n \sigma_i^2(A) u_i\right)d
d\sigma_n.
\]
\end{corollary}
\label{reduction}
\begin{proof}
Lemma \ref{invar} tells us that the inner integral does not depend on $B \in
\mathcal{A},$ and Lemma \ref{singinv} tells us how to evaluate that integral
in terms of singular values of $A.$
\end{proof}
Corollary \ref{reduction} tells us that Question \ref{wilk} will be answered
once we succeed in showing Theorem \ref{mainineq} (for diagonal $A.$)

\begin{theorem-definition}[Gershgorin's Theorem]
Let $A \in M_n(\cx).$ Let 
\begin{gather}
r_i = \sum_{j\neq i} |a_{ij}|,\\
s_i = \sum_{j\neq i} |a_{ji}|,
\end{gather} and let the row- and column- \emph{Gershgorin disks}
(respectively) be:
\begin{gather}
\mathcal{R}_i = \{z\in \cx \quad | \quad |z - a_{ii}| \leq r_i\},\\
\mathcal{S}_i = \{z\in \cx \quad | \quad |z - a_{ii}| \leq s_i\}.
\end{gather}
Let now $\lambda$ be an eigenvalue of $A.$ Then there exists a $k$ and an $l$
such that $\lambda \in \mathcal{R}_k$ and $\lambda \in \mathcal{S}_l.$
\end{theorem-definition}
\begin{proof}
See, eg, \cite{bm94} or \cite{hj90}.
\end{proof}

\section{Proof of the Theorem \ref{mainineq}}
\label{mainproof}

To prove Theorem \ref{mainineq}, 
by Lemma \ref{singinv}, it will be sufficient to prove
\begin{theorem}
\label{redineq}
\[
\mathfrak{I}(a_1, \dotsc, a_n) = \int_\essn \log\left(\sum_{i=1}^n a_i
u_i^2\right) d \sigma_n \geq 0,
\]
whenever $a_1, \dotsc a_n > 0,$ and $\prod_{i=1}^n a_i = 1.$ Furthermore,
equality holds only when $a_1 = \cdots = a_n = 1.$
\end{theorem}

\begin{proof}
By the (weighted version of) the arithmetic-geometric mean inequality (see
\cite{hlp88}),  
\[
\sum_{i=1}^n a_i x_i^2 \geq \prod_{i=1}^n a_i^{x_i^2},
\]
with equality if and only if all of the $a_i$ are equal.
(the above uses the fact that $\sum_{i=1}^n x_i^2 = 1.$)
Taking logs, we see that
\[
\log\left(\sum_{a_i x_i^2}\right) \geq \sum_{i=1}^n x_i^2 \log(a_i),
\]
with equality if and only if all the $a_i$ are equal.
And integrating, we see that 
\begin{multline*}
\mathfrak{I}(a_1, \dotsc, a_n) \geq \int_\essn \sum_{i=1}^n x_i^2 \log(a_i)
= \\
\sum_{i=1}^n \log(a_i) \int_\essn x_i^2 d\sigma_n = 
 n \log(a_i) \int_\essn x_1^2 d\sigma_n = 0,
\end{multline*}
with equality if and only if $a_i = 1$ for all $i.$
\end{proof}

\section{A lower bound on average spectral radius}
\label{lbound}
In this section we consider the following average:
\[
\mathfrak{A}(A) = \fint_{O_n} \log \rho(A X) d X,
\]
where $\rho(M)$ denotes the spectral radius of the matrix $M.$
The main theorem of the section will be the following estimate:
\begin{theorem}
\label{avthm}
Let $A \in SL(n, \mathbb{R}),$  let $n \geq 3,$ and let $\sigma_1(A)$ be the
largest singular value of $A.$ 
Then there exists a constant $C_n > 0,$ such that 
\begin{equation}
\label{aveqn}
\mathfrak{A}(A) \geq C_n \log \sigma_1(A).
\end{equation}
\end{theorem}
\begin{remark} Theorem \ref{avthm} is \emph{false} for $n=2$ and trivially
  true for $n=1$ ($C_1 = 1$).
\end{remark}
To prove Theorem \ref{avthm} we first make a couple of observations. First, 
\begin{lemma}
\label{singsuf}
Let $A \in M_n,$ let $\Sigma(A)$ be the singular values of $A,$ and let
$\rho(X)$ be the spectral radius of the matrix $X.$
Then
\[
\fint_{O_n} f(\rho(AY)) d Y = \fint_{O_n} f(\rho\left(\diag\left(\Sigma
A\right)\right) d Y.
\]
\end{lemma}
\begin{proof}
Form the singular value decomposition of $A,$ that is, write 
$A = P_1 \diag\left(\Sigma(A)\right) P_2,$ where $P_1, P_2 \in O_n.$ Then, 
\[
\rho\left(A Y\right) = \rho\left(\diag\left(\sigma(A)\right) P_2 Y P_1\right).
\]
The statement of the Lemma follows immediately.
\end{proof}
Lemma \ref{singsuf} tells us that in the statement of Theorem \ref{avthm}, it
is enough to consider \emph{diagonal} matrices $A$ with \emph{positive}
entries. 
In the rest of this 
section we will labor under these assumptions. The strategy will be to first
demonstrate a bound in any closed set of diagonal matrices not containing the
identity matrix $I_n$ -- this will be done in Section \ref{away} (particularly
Corollary \ref{boundaway}--  and then show a bound in some small neighborhood of the
identity -- this will be the content of Section \ref{hiho}. Our measure of the
distance from the identity will be the quantity 
\[
L(A) = \max_{i=1}^n |\log(\sigma_i(A)|.
\]
the top singular value $\sigma_1(A).$ From now on,
$A=\diag(d_1, \dots, d_n),$ with $d_1 \geq d_2 \geq \dots d_n > 0.$ We will
assume that $\prod_{i=1}^n d_i =1$ when 
necessary -- much of the time we will only use the fact that $d_1 > 1.$

\section{Away from identity}
\label{away}
\begin{theorem}
Let $A = \diag(d_1, \dots, d_n),$ with $d_1 \geq d_2, \dotsc, \geq d_n \geq 0,$ and
  let $\|M\|_1 \leq 1.$ Let $A_t = A(I+tM).$ Then,
\[
d_1(1+ 2 t) \geq \rho(A_t) \geq d_1(1-2 n t).
\]
\end{theorem}
\begin{proof}
Since $(A_t)_{ij} = d_i (\delta_{ij} + tM_{ij}),$ it follows that 
\begin{equation}
\label{g1}
d_i(1+ t) \geq |(A_t)_{ii}| \geq d_i(1-t),
\end{equation}
while
\begin{equation}
\label{g2}
|r_i(A_t)| \leq d_i t.
\end{equation}
Let \[\mathcal{G}(A_t) = \bigcup_{i=1}^n R_i(A_t).\]
By a simple continuity argument (see, for example, \cite{bm94} for a
similar argument), the perturbation of the eigenvalue $d_1$ of $A$ stays in
the connected component of the Gershgorin disk $\mathcal{R_1}(A_t).$ Since by
the above inequalities Eq.~\eqref{g1} and Eq.~\eqref{g2}, 
\begin{gather*}
\min_{z\in \mathcal{R}_i(A_t)}  \geq d_i(1-2  t),\\
\max_{z\in \mathcal{R}_i(A_t)}  \leq d_i(1-2  t),\\
\end{gather*}
we see that after time $t$ the connected component of $\mathcal{R}_1(A_t)$ can
  only stay between the advertised bounds.
\end{proof}

\begin{theorem}
\label{farest}
Let $A = \diag(d1, \dots, d_n).$ 
For any $\epsilon > 0,$ there is a constant $c_\epsilon > 0,$ such that 
\[
\mathfrak{A}(A) \geq c_\epsilon(d_1 - \epsilon).
\]
\end{theorem}
\begin{proof}
First, let $t_0$ be so small that $\log(1-2 n t_0) > 1- \epsilon.$
Now let 
\[
S_0 = \{x \in O_n \quad | \quad \| x - I\|_1 < t_0\},
\]
and let 
\[
c_\epsilon = \dfrac{\vol(S_0)}{\vol(O_n)}.
\]
We see that 
\[
\mathfrak{A}(A) \geq 
\dfrac{\int_{S_0}
\log \rho(A X) d X}{\vol O_n} \geq c_\epsilon(d_1 - \epsilon).
\]
\end{proof}
\begin{corollary}
\label{boundaway}
Let $\mathcal{D}$ be the set of all diagonal matrices with elements $d_1 \geq
d_2 \geq \dots \geq d_n > 0,$ and let $\mathcal{D}_t$ be the set of those
matrices $M \in \mathcal{D}$ with $d_1(M) \geq t > 1.$ Then there is a
constant $C_t > 0$ such that 
\[
\dfrac{\mathfrak{A}(A)}{\log(d_1(A))} > C_t,
\] for all $A \in \mathcal{D}_t.$
\end{corollary}
\begin{proof}
Let $s = \log(t)/2.$ Then, by Theorem \ref{farest}, 
\[
\mathfrak{A}(A) \geq c_s (\log(d_1(A) - s) \geq \dfrac{c_s \log d_1(A)}2,
\]
where the last inequality holds for all $A \in \mathcal{D}_t.$ Now set $C_t =
c_s/2.$
\end{proof}

\section{Close to the Identity}
\label{hiho}

In this section we shall prove the following:
\begin{theorem}
\label{boundclose}
Let $n \geq 3.$ 
Let $A = \diag(\exp(d_1 t), \dotsc, \exp(d_n t)),$ with
$d_1 \geq \dots \geq  d_n,$ with  $\max(|d_1|, |d_n|) = 1.$ and 
\[
\sum_{i=1}^n d_i = 0.
\]
Then there exist constants $\mathcal{C}_n$  and $\epsilon_n,$ such that for $t < \epsilon_n,$ 
\begin{equation}
\label{closebound}
\mathfrak{A}(A) \geq \mathcal{C}_n t d_1
\end{equation}
for $0<t \leq \epsilon_n.$
\end{theorem}

\begin{remark}
The argument which follows can be made to give quite explicit lower bounds for
both $\epsilon_n$ and $\mathcal{C}_n.$ However, this does not seem useful,
since it is quite clear that the bounds will not be close to the truth. It
seems plausible that a combination of the methods of this note with a detailed
understanding of the structure of the orthogonal group will give a tight
result, but we shall not attempt to do this here. It should, on the other
hand, be noted that not using any structural property of the orthogonal group
means that the arguments in this note work just as well for cosets of unitary
groups. 
\end{remark}

To prove Theorem \ref{boundclose} we shall first need the following trivial
but crucial observation:
\begin{observation}
\label{localized}
Let $A \in SL(n, \reals), $ and let $X \in O_n(\reals).$ Then $\rho(A X) \geq
1.$
\end{observation}
\begin{proof}
$\det A X = 1.$
\end{proof}

Now, recall that 
\[\mathcal{A}(A) = \fint_{O_n} \log \rho (AX) d x.\]
By the observation above, it follows that the integrand is everywhere
nonnegative, and so to prove a lower bound such as that of Theorem
\ref{boundclose} it will be enough to show the inequality \eqref{closebound}
with $\mathcal{A}(A)$ replaced by 
\[
\mathcal{A}_S(A) = \int_{S\subseteq O_n} \log \rho(AX) d x.
\]
Let us pick $S$ to be a small neighborhood around a matrix $X_0
\in O_n.$

Notice that all eigenvalues of $X \in O_n$ have absolute value $1.$ We will chose
$X_0$ to have a \emph{simple} eigenvalue $\lambda_0,$ such that
\begin{equation}
\label{localineq}
\dfrac{d \log |\lambda_0(\diag(\exp(t d_1), \dotsc, \exp(t d_n))X_0)|}{d t}
\geq C |d_1|,
\end{equation}
for some constant $C.$ By analyticity of $\lambda_0$ (Theorem
\ref{analyticity}),  inequality \eqref{localineq} suffices to show that 
\begin{equation}
\label{closeboundloc}
\mathfrak{A}(A) \geq \mathcal{C}_n t d_1,
\end{equation}
for $t$ sufficiently small, and by the Observation \ref{localized} this will
suffice to show Theorem \ref{boundclose}.
\subsection{Proof of inequality \eqref{localineq}}
First, we need to find a suitable orthogonal matrix $X_0.$ There will be two
cases, depending on the parity of $n.$

First case is that of odd $n.$ In such a case, we pick $X_0$ to be a rotation
by $180$ degrees around the axis $e_1 = (1, 0, \dotsc, 0),$ and our
$\lambda_0(X_0) = 1.$ Its eigenvector is $v_{\lambda_0}= e_1.$ 
By Lemma \ref{expvar},
\[
\dfrac{\lambda(x)}{dx} = \lambda(x)  d_1,
\]
and we are done.

Second case is that of $n$ even. In that case, we pick $X_0$ to fix the span
of $e_3, \dots, e_n$ and to rotate the span of $e_1$ and $e_2$ by $90$
degrees.
In this case, $\lambda_0(X_0) = i,$ and the eigenvector of $\lambda_0$ is
\[
v_{\lambda_0}=\dfrac{1}{\sqrt{2}}(1, i, 0, \dotsc, 0).
\]

By Lemma \ref{expvar},

\[
\dfrac{\lambda(x)}{dx} = \lambda(x) \frac{1}{\sqrt{2}}( d_1+d_2).
\]
If $d_2 > 0,$ we are done, since $d_1 + d_2 > d_1.$ Since
\[
\sum_{i=1}^n d_i = 0,
\]
it follows that 
\[
d_1 = -\sum_{i=2}^n d_i = \sum_{i=2}^n |d_i|.
\]
Since $|d_n| \geq |d_{n-1} \geq |d_2|,$ it follows that $d_1 \geq (n-1)
|d_2|,$ and so 
\[
d_1 + d_2 \geq \dfrac{n-2}{n-1} d_1,
\]
so as long as $n > 2$ we are done. Notice that the argument breaks down when
$n=2,$ and in fact, in that case the result is false.
\section{What about dimension 2?}
\label{dim2}
The function $f(d_1, \dotsc, d_n) = \mathfrak{A}(\diag(\exp(d_1), \dotsc,
\exp(d_n)))$  restricted to the set 
\[
\sum_{i=1}^nd_1\dots d_n = 0
\]
always has a minimum at the origin. In dimension $2,$ however, the minimum is
an actual critical point. This is not hard to see using Lemma \ref{expvar},
together with the observation that the group $O(2)$ is abelian, and thus every
$x\in O(2)$ has $(1, i)$ and $(i, 1)$ for an orthonormal basis of
eigenvectors. Thus, at $d_1 = 0,$ 
\[
\dfrac{d\log(\rho(\diag(\exp(d_1), \exp(-d_1)) X)}{d d_1} = 0
\]
for \emph{every} $X \in O(2).$
In dimension $n>2$ the function $\mathfrak{A}$ is not smooth at the origin.

As a matter of fact, in dimension 2 many of the computations can be made much
more explicit. Indeed, a general element $R_{\theta}$ of $SO(2)$ has the form:
\[
R_\theta = \begin{pmatrix}
\cos(\theta) & \sin(\theta)\\
-\sin(\theta) & \cos(\theta)
\end{pmatrix},
\]
and so if $A = \diag(a, 1/a),$ then, denoting $t = (a+1/a)/2,$
\[
\tr A R_{\theta} = 2 t \cos(\theta),
\]
and so 
\[
\rho(A R_{\theta} = \max\left(|t\cos(\theta)+\sqrt{t^2 \cos^2(\theta) -1}|,
|t\cos(\theta)-\sqrt{t^2 \cos^2(\theta) -1}|\right).
\]
A simple computation shows:

\begin{lemma}
\label{support}
Whenn $|\cos(\theta)|< 1/t,$ then $\rho(A
R_{\theta}) = 1.$ (in other words, the matrix $A R_{\theta}$ is
\emph{elliptic}. 
\end{lemma}

We have the following corollary (which answers a question of \cite{lssw03}):
\begin{corollary}
\label{sup2}
Let $\mu$ be a measure on $SO(2)$ such that the support of $\mu$ omits some
interval. Then, there exists a matrix $A \in SL(2, R),$ such that
\[
\int_{SO(2)} \log \rho(A R_{\theta}) d\mu = 0.
\]
\end{corollary}
\begin{proof}
Suppose that the support of $\mu$ omits $[a, b].$ Let $c = (a+b)/2$ ($a, b, c$
are all modulo $2\pi.$) Let $s = |(a-b)/2|,$ let $\alpha = 1/\cos(s).$ and let
$\beta$ be positive, and such that $\beta + 1/\beta < 2\alpha$ (in other
words, $\beta < \exp(\acosh(\alpha)).$)
Let $A = \diag(\beta, 1/\beta) R_{-c}.$ Then Lemma \ref{support} shows that
$A$ is the matrix whose existence is asserted in the corollary.
\end{proof}

\section{Some perturbation theory}
\label{perturb}
In this section we will recall some formulas of perturbation theory (an
exhaustive treatment can be found in the classic \cite{kato95}), and derive
some estimates needed in the rest of this paper.

\begin{theorem}
\label{analyticity}
Let $T$ be an $n \times n$ matrix (thought of as a point in $\cx^{n^2},$ and
let $\lambda$ be a \emph{simple} eigenvalue of $T.$ Then there is a
neighborhood of $T$ where $\lambda$ is a holomorphic function of $T.$ In
particular, $\lambda$ is $C^\infty$ in a neighborhood of $T.$
\end{theorem}
\begin{proof}
This follows from generalities on analyticity of (simple) roots of
polynomials as functions of their coefficients. See \cite{kato95} for some
further details.
\end{proof}

\begin{theorem}
\label{katoform}
Let $\lambda$ be a simple eigenvalue of $T(0).$ Then 
\[
\dfrac{d \lambda}{d x}\left|_{x=0}\right. =
\tr \left(\dfrac{T(x)}{dx}\left|_{x=0}\right. P_\lambda\right),
\]
where $P_\lambda$ is the orthogonal projection onto the eigenspace of
$\lambda.$
\end{theorem}
\begin{proof}
See \cite[p.~79]{kato95}.
\end{proof}

\begin{lemma}
\label{expvar}
Let $T(x) = \diag(\exp(d_1 x), \dotsc, \exp(d_n x)) T,$ and let let
$\lambda(x)$ be a simple eigenvalue of $T(x).$ Let the unit eigenvector of
$\lambda$ be $v_\lambda = (v_1, \dotsc, v_n).$ Then,
\[
\dfrac{\lambda(x)}{dx} = \lambda(x)  \sum_{i=1}^n d_i |v_i|^2.
\]
In particular, the logarithmic derivative of $\lambda(x)$ is real.
\end{lemma}

\begin{proof}
Use Theorem \ref{katoform} and the observation that \[T(x) P_{\lambda(x)} = 
\lambda(x) P_{\lambda(x)},\]
to write
\[
\dfrac{\lambda(x)}{dx} =\tr \left(\diag(d_1, \dotsc, d_n)
P_{\lambda(x)}\right).\] 
Now, note that $(P_{\lambda(x)})_{ij} = v_i \overline{v_j}.$
\end{proof}

\bibliographystyle{plain}
\bibliography{amie}
\end{document}